\documentclass[12pt]{article}
\usepackage{amsmath,amssymb}

\newcommand{\prf}{\noindent{\bf Proof.}\ }
\newcommand{\Tbar}{\ensuremath{\overline{\mathcal{T}}}}

\newcommand{\caT}{\mathcal{T}}
\newcommand{\caC}{\mathcal{C}}
\newcommand{\geod}{\widehat{p_1T^np_2}}

\newtheorem{definition}{Definition}[section]
\newtheorem{lemma}[definition]{Lemma}
\newtheorem{theorem}[definition]{Theorem}

\newtheorem{proposition}[definition]{Proposition}

\begin{document}

\title{Topological dynamics of the Weil-Petersson geodesic flow\footnote{2000 Mathematics Subject Classification Primary: 37D40, 32G15; Secondary: 53D25.}}         
\author{Mark Pollicott, Howard Weiss and Scott A. Wolpert}       
\date{June 1, 2009}          
\maketitle

\begin{abstract}
We prove topological transitivity for the Weil-Petersson geodesic flow for real two-dimensional moduli spaces of hyperbolic structures.  Our proof follows a new approach that combines the density of singular unit tangent vectors, the geometry of cusps and convexity properties of negative curvature.  We also show that the  Weil-Petersson geodesic flow has: horseshoes, invariant sets with positive topological entropy, and that there are infinitely many hyperbolic closed geodesics, whose number grows exponentially in length.  Furthermore, we note that the volume entropy is infinite.
\end{abstract}
\section{Introduction}

The moduli space $\mathcal M_{g,n}$ is the space of hyperbolic metrics for a surface of genus $g$ with $n$ punctures. There are several interesting metrics on $\mathcal M_{g,n}$, and during the past few years there has been intense activity on studying the geometry and relationships between these metrics \cite{LSY2}. There has also been significant activity on investigating the fine dynamics of the geodesic flow for the Teichm\"{u}ller metric \cite{AGY, AF, ABEM}.  The Teichm\"{u}ller metric is a complete Finsler metric that describes a straight space, geodesics connecting points are unique and maximal geodesics are bi infinite.  
Minsky has provided a simple model for the Teichm\"{u}ller metric.  
For a Riemann surface with a collection of non peripheral, disjoint, distinct free homotopy classes $\mathcal F$, parameterize Riemann surfaces by: the (product) Teichm\"{u}ller space for the $\mathcal F$ complement and the Fenchel-Nielsen parameters for $\mathcal F$.  Introduce a comparison metric - the supremum of the Teichm\"{u}ller metrics for the components of the complement, and formal hyperbolic metrics for the Fenchel-Nielsen length-angle planes.  The comparison metric approximates within an additive constant in the region of Teichm\"{u}ller space, where $F$ represents the short hyperbolic geodesics \cite{Minel}.    

We undertake an investigation of the topological dynamics of the geodesic flow for the Weil-Petersson (WP) metric.  Our results are for real two-dimensional moduli spaces, i.e., the moduli space $\mathcal M_{1,1}$ for once punctured tori and the moduli space $\mathcal M_{0,4}$ for spheres with four punctures. These spaces are completed by adding points, {\em cusps}, for the degenerate hyperbolic structures, with the property that a dense set of singular unit tangent vectors are initial to geodesics that end at a cusp in finite time \cite{Brkwpvs}. Although non-compact, these moduli spaces have finite area and finite diameter.  Curvature is negative, bounded away from zero and is unbounded in each cusp \cite{Wlcurvnt}.  It follows with these properties that the geodesic flow (GF) is an incomplete uniformly hyperbolic flow, and thus needs to be studied as a non-uniformly hyperbolic dynamical system with singularities.

Our main result is topological transitivity for the WP GF for two-dimensional moduli spaces.  Our proof is based on a new approach.  Earlier proofs of topological transitivity for geodesic flows or systems with singularities use some form of coding, by introducing an ideal boundary or a Markov partition. Given the density of singular tangent vectors, there appears to be no complete notion of an ideal boundary \cite{Brkwpvs}, and the construction of a Markov partition would be quite delicate. The geometry of a cusp combined with the $CAT(-\epsilon),\,\epsilon>0,$ geometry provide a substitute for coding. We develop a {\em shadowing lemma}, Proposition \ref{distbd}, to approximate piecewise geodesics by geodesic chords.  We find that the WP GF on $\mathcal M_{1, 1}$ and $\mathcal M_{0, 4}$ has: horseshoes, subsets with positive topological entropy, and that there are infinitely many hyperbolic closed geodesics, whose number grows exponentially in length.  Furthermore, we find that the volume entropy is infinite.

We start considerations with the density in the unit tangent bundle of $\mathcal M_{g,n}$ of geodesics connecting cusps of $\mathcal T_{g,n}$.  We next use that in dimension two a pair of geodesic segments emanating form a cusp (modulo Dehn twists) is approximated by geodesic segments connecting a point and Dehn twist translates of a point.  A piecewise geodesic $\mathcal C$ of controlled small exterior angle is then constructed from geodesics connecting cusps and approximating geodesic segments.  A limit of chords of $\mathcal C$ is a bi infinite geodesic asymptotic to $\mathcal C$.  We use the shadowing lemma to control the chords and to establish the asymptotic behavior.  

Brock, Masur and Minsky establish topological transitivity for all moduli spaces $\mathcal M_{g,n}$ \cite{BMM}.  They start by defining an {\em ending lamination} for an infinite WP geodesic ray.   They find for the full measure set of bi infinite WP geodesics, that bi recur to some compact set of $\mathcal M_{g,n}$, that the ending laminations characterize the ending asymptote classes.  Given a geodesic connecting cusps, they approximate by bi recurrent geodesics.  For a bi recurrent geodesic they use the resulting ending laminations as data for specification of a pseudo Anosov mapping class with axes approximating the geodesic.  They further use the dynamics of pseudo Anosov elements and compositions to approximate a sequence of geodesics.  The main result follows.  They also show that axes of pseudo Anosov elements are dense in the unit tangent bundle of $\mathcal M_{g,n}$.    

The present and Brock, Masur and Minsky approaches begin with the density of geodesics connecting cusps.  The approaches differ in that the present approach is based on approximating piecewise geodesics and a shadowing argument, whereas the later approach develops a partial boundary theory to construct pseudo Anosov elements with controlled axes.

Our approach also applies to complete constant negatively curved surfaces with cusps. We believe that the approach should also have applications to certain systems with singularities, including billiards.

\section{Preliminaries}
We present a unified discussion for the basic metric geometry of the Weil-Petersson metric.  The treatment is based on earlier and recent works of a collection of authors including \cite{Abbook, Bersdeg, Brkwpvs,DW2,Zh2,MW,Wlcomp,Wlbhv,Wlqk}.  Once punctured tori and four punctured spheres are related in the tower of coverings.  The Teichm\"{u}ller spaces are isomorphic and isometric.  We study once punctured tori.

Let $\caT$ be the Teichm\"{u}ller space for marked once punctured tori.  Points of $\caT$ are equivalence classes $\{(R,ds^2,f)\}$ of marked complete hyperbolic structures $R$  with reference homeomorphisms $f:F\rightarrow R$ from a base surface $F$.  A self homeomorphism of $F$ induces a mapping of $\caT$ and the mapping class group $MCG$ of $F$ acts properly discontinuously on $\caT$ with quotient the moduli space $\mathcal M$ of Riemann surfaces.  The Teichm\"{u}ller space $\caT$ is a complex one-dimensional manifold with cotangent space at $R$ being $Q(R)$, the space of holomorphic quadratic differentials with at most simple poles at punctures.  

The Teichm\"{u}ller-Kobayashi co metric for $\caT$ is defined as
\[
\|\varphi\|_T\,=\,\int_R|\varphi| \quad\mbox{for}\quad \varphi\in Q(R).
\]
For complex one-dimensional Teichm\"{u}ller spaces the Teichm\"{u}ller metric coincides with the hyperbolic metric.  The Weil-Petersson (WP) co metric is defined as
\[
g_{WP}(\varphi,\psi)\,=\,\int_R\varphi\overline{\psi}\, (ds^2)^{-1}\quad\mbox{for}\quad \varphi,\psi\in Q(R)
\]
and $ds^2$ the hyperbolic metric of $R$.  The Teichm\"{u}ller and WP metrics are $MCG$ invariant.

A non trivial, non peripheral, simple free homotopy class on $F$, now called an {\it admissible} free homotopy class, determines a geodesic $\gamma$ on $R$ and a hyperbolic pants decomposition.  The length $\ell_{\gamma}$ for the joined pants boundaries and the offset, or {\it twist} $\tau_{\gamma}$, for adjoining the pants boundaries combine to provide the Fenchel-Nielsen (FN) parameters $(\ell_{\gamma},\tau_{\gamma})$ valued in $\mathbb R_+\times\mathbb R$.  The FN parameters provide global real analytic coordinates with each admissible free homotopy class and choice of origin for the twist determining a global coordinate.  The Dehn twist $T$ about $\gamma$, an element of $MCG$, acts by $T\colon(\ell_{\gamma},\tau_{\gamma})\to(\ell_{\gamma},\tau_{\gamma}+\ell_{\gamma})$.    

A bordification (a partial compactification) $\Tbar$, the {\it augmented Teichm\"{u}ller space}, is introduced by extending the range of parameters.  For an $\ell_{\gamma}$ equal to zero, the twist $\tau_{\gamma}$ is not defined and in place of the geodesic $\gamma$ on $R$ there appears a further pair of cusps (a formal node).  The extended FN parameters describe marked (possibly) noded Riemann surfaces.  The points of $\Tbar-\caT$ are the {\it cusps of Teichm\"{u}ller space} and are in one-one correspondence with the admissible free homotopy classes.  The neighborhoods of $\{\ell_{\gamma}=0\}$ in $\Tbar$ are given as $\{\ell_{\gamma}<\epsilon\}$, $\epsilon>0$; $\Tbar$ is non locally compact, since neighborhoods of $\{\ell_{\gamma}=0\}$ are stabilized by the Dehn twist about $\gamma$; $\caT\subset\Tbar$ is a convex subset. Each pair of distinct cusps of $\Tbar$ is connected by a unique finite length unit-speed WP geodesic, now called a {\it singular geodesic.}  There is a positive lower bound for the length of singular geodesics.  Tangents to singular geodesics at points of $\caT$ are called {\it singular tangents.}

The augmented Teichm\"{u}ller space $\Tbar$ is $CAT(0)$  and in particular is a unique geodesic length space, \cite{BH}.  Triangles in $\Tbar$ are approximated by triangles in $\caT$ and for once punctured tori the WP curvature is bounded above by a negative constant \cite{Wlcurvnt}.  In particular for once punctured tori,  the curvature of the augmented Teichm\"{u}ller space is bounded away from zero, and thus the space is $CAT(-\epsilon)$ for  some $\epsilon>0$.     

\begin{theorem}
\label{bkrnd}
The augmented Teichm\"{u}ller space for once punctured tori $\Tbar$ is the WP completion of $\caT$ and is a $CAT(-\epsilon)$ space for some $\epsilon>0$. The WP curvature is bounded above by a negative constant and tends to negative infinity at the cusps of Teichm\"{u}ller space.  The singular geodesics have tangents dense in the unit tangent bundle of Teichm\"{u}ller space.
\end{theorem}

For once punctured tori the upper half plane $\mathbb H$ and elliptic modular group $PSL(2;\mathbb Z)$ are the models for $\caT$ and $MCG$.  The augmentation $\Tbar$ is $\mathbb H\cup \mathbb Q$ provided with the horocycle topology.  The rational numbers are the cusps of $\caT$.  The density for singular WP geodesics \cite[Coro. 18]{Wlcomp} is a counterpart of the density for the hyperbolic metric of geodesics with rational endpoints.

We are interested in the geometry of the cusps of Teichm\"{u}ller space.  The WP metric has the cusp expansion
\[
g_{WP}\,=\,2\pi((d\ell_{\gamma}^{1/2})^2 \,+\,(d\ell_{\gamma}^{1/2}\circ J)^2)\,+\, O(\ell_{\gamma}^3\,g_{WP})
\]
for $J$ the almost complex structure, and for bounded length $\ell_{\gamma}$ for the curve corresponding to the cusp.  We are most interested in geodesics that approach a cusp.  A unit tangent at a point of $\caT$ is {\it special} for a cusp provided it is tangent to the geodesic ending at the cusp (a tangent $v$ is singular provided $v$ and $-v$ are special).  For a given cusp, at each point of $\caT$ there is a unique associated special unit tangent, the initial tangent for the geodesic connecting to the cusp.  We consider the spiraling of geodesics to a cusp as follows.

\begin{proposition}
For a given admissible free homotopy class, let $T$ be the Dehn twist and $v_1,v_2$ special tangents with base points $p_1,p_2$.  Given $\epsilon>0$ there exists a positive integer $n_0$,  such that for $n\ge n_0$ the geodesic $\geod$ connecting $p_1$ to $T^np_2$ has initial and terminal tangents respectively within $\epsilon$ of $v_1$ and $-T^nv_2$.
\end{proposition} 
\prf  The behavior of $\ell_{\gamma}$ on the geodesic $\geod$ is basic.  The geodesic-length function $\ell_{\gamma}$ is $T$-invariant and convex on $\geod$.  The function $\ell_{\gamma}$ is bounded by its values at $p_1$ and $p_2$.  The Dehn twist $T$ acts in Fenchel-Nielsen coordinates by $(\ell_{\gamma},\tau_{\gamma})\to(\ell_{\gamma},\tau_{\gamma}+\ell_{\gamma})$.  
We let $c_1=\inf \ell_{\gamma}(\geod)$ ($c_1$ could be zero) and 
$c_2=\max \{\ell_{\gamma}(p_1),\ell_{\gamma}(p_2)\}$.  We further let $L_0$ be the infimum of WP length for geodesics connecting points $q_1,q_2$ satisfying $c_1\le \ell_{\gamma}(q_1), \ell_{\gamma}(q_2) \le c_2$ and $\tau_{\gamma}(q_2)=\tau_{\gamma}(q_1)+1$.  By compactness of the set of connecting geodesics, $L_0$ is positive provided $c_1$ is positive.  A segment of $\geod$ with $a\le \tau_{\gamma}\le a+1$ satisfies the stated conditions and so has length at least $L_0$.  It follows that the length satisfies $L(\geod)\ge (n-n_0)L_0$, for a suitable $n_0$.    
  At the same time, if we write $p$ for the cusp $\{\ell_{\gamma}=0\}$, then $L(\geod)\le d(p_1,p)+d(p,p_2)$ since $\widehat{p_1p}$, $\widehat{pT^np_2}$ combine to connect $p_1$ to $p_2$.  We have for $c_1\le\ell_{\gamma}\le c_2$ on $\geod$ that $(n-n_0)L_0\le d(p_1,p)+d(p,p_2)$. It now follows for $n$ tending to infinity that $\min \ell_{\gamma}$ on $\geod$ tends to zero.  

We write $q_n$ for the minimum point for $\ell_{\gamma}$.  By definition of the topology of $\Tbar$, the points $q_n$ and $T^{-n}q_n$ limit to the cusp $p$.  For a $CAT(0)$ space the distance between geodesic segments is convex and achieves its maximum at endpoints.  It follows that  $\widehat{p_1q_n}$ tends uniformly to $\widehat{p_1p}$ and $\widehat{T^{-n}q_np_2}$ tends uniformly to $\widehat{pp_2}$.  On $\caT$ local convergence of geodesics implies $C^1$-convergence and the desired convergence of tangents.  The proof is complete.

The counterpart result for the upper half plane is valid.  For the cusp at infinity, consider a translation $T$ and vertical vectors at the basepoints $p_1,p_2$.  The vectors are approximated by the initial and terminal tangents of the geodesics $\geod$.   

Generalizing the proposition for higher-dimensional Teichm\"{u}ller spaces is an open question.  Approximating a pair of special tangents by the initial and terminal tangents of connecting geodesics is only possible for a proper subset of tangent pairs.  A special geodesic has a non zero terminal tangent in the Alexandrov tangent cone of the terminal point in $\Tbar$, \cite{Wlbhv}.  In Example 4.19 of \cite{Wlbhv}, a necessary condition is provided for a pair of terminal tangents to arise from a twisting limit of connecting geodesics. The expectation is that a necessary and sufficient condition for approximating a pair of special tangents is for the terminal tangents to satisfy the condition of Example 4.19.  

{\bf Twisting of geodesics.}  We  describe a piecewise-geodesic modification of a pair of geodesics $\gamma_1,\gamma_2$ emanating from a cusp of Teichm\"{u}ller space.  Given distance parameter $\delta>0$, we will replace the segments of $\gamma_1,\gamma_2$ within $\delta$ of the cusp to obtain a concatenation of three geodesic segments with small exterior angles at the two concatenation points.  Let $p_j$ be the point of $\gamma_j$ at distance $\delta$ from the cusp.  From Proposition 2.2,  for the twisting exponent $n$ sufficiently large, the exterior angles at $p_1$ (and at $T^np_2$) between $\gamma_1$ and $\widehat{p_1T^np_2}$ (at $T^np_2$ between $T^n\gamma_2$ and $\widehat{p_1T^np_2}$) are sufficiently small.  {\bf Conclusion:} for $n$ sufficiently large, then the concatenation of three segments - the segment of $\gamma_1$ outside the $\delta$-ball {\em concatenate} $\widehat{p_1T^np_2}$ {\em concatenate} the segment of $T^n\gamma_2$ outside the $\delta_1$-ball - is piecewise-geodesic with sufficiently small exterior angles.  Furthermore, given a smaller neighborhood of the cusp, for $n$ sufficiently large, the small exterior angles provide that outside the smaller neighborhood the segments of $\widehat{p_1T^np_2}$ are $C^1$ sufficiently close to their corresponding segments of $\gamma_1, \gamma_2$.  In particular for $n$ sufficiently large, the concatenation is sufficiently close to $\gamma_1\cup T^n\gamma_2$.      

\section{Piecewise-geodesics and chordal limits}
We construct infinite WP geodesics dense in the unit tangent bundle of the moduli space. The prescription begins with a  sequence of WP singular geodesics on $\caT$ with projection dense in the unit tangent bundle.  We then apply twisting to construct a sequence of piecewise geodesics with small exterior angles, and consider a limit of chords for points tending to infinity.  We preclude degeneration of the chords and using a shadowing lemma, we establish the existence of a suitable limit, provided the sequence of exterior angles has small norm in the sequence space $\ell^1$.  The result is presented in Theorem \ref{dense}.  We begin consideration with observations about piecewise-geodesics with small exterior angles.  

We begin with the discussion of chords for concatenations of geodesics.  The first consideration is the distance to a closed geodesic segment.   The Teichm\"{u}ller space $\caT$ is a convex subset of a $CAT(-\epsilon), \epsilon>0$ space and is a Riemannian surface.  Accordingly for a geodesic segment $\gamma$, closed in $\caT$, we introduce the nearest point projection $\Pi$ to $\gamma$, \cite[Chap. II.2]{BH}.  The fiber of $\Pi$ for an interior point of $\gamma$ is the complete geodesic orthogonal to $\gamma$ at the basepoint of the projection.  The fiber of $\Pi$ for an endpoint is the domain bounded by the geodesic orthogonal to the endpoint.  We recall the formula for the derivative of the distance to $\gamma$ along a smooth curve $\sigma$.  At a point $p$ of $\sigma$, consider the angle $\theta$ between $\sigma$ and the away pointing fiber of the projection.  The derivative of the distance $d_{\gamma}$ to $\gamma$ along the unit-speed parameterized curve $\sigma$ is 
$\cos \theta$.  

We consider a (possibly infinite) concatenation $\mathcal C$ of geodesic segments with vertices in $\caT$, and write $ea(\caC)$ for the finite total variation of exterior angles at vertices.  We first observe for $\gamma$ a geodesic chord between points of $\caC$ the bound $\max |d_{\gamma}'|\le ea(\caC)$, for the derivative of distance to $\gamma$.  The argument is as follows.  The distance $d_{\gamma}$ is convex along the segments of $\caC$ with jump discontinuities for $d_{\gamma}'$ at vertices.  In particular the graph of $d_{\gamma}'$ consists of increasing intervals and jump discontinuities.  Since the cosine is $1$-Lipschitz the total negative jump discontinuity of $d_{\gamma}'$ is at most $ea(\caC)$.  
Consider a segment $\widehat{pq}$ of $\gamma$ intersecting the concatenation $\caC$ only at endpoints.  The function $d_{\gamma}$ is positive on $\widehat{pq}$ and so has a non negative initial derivative and non positive final derivative.  It follows for $\widehat{pq}$, the absolute total jump discontinuity of $d_{\gamma}'$ bounds the sum of the initial value of $d_{\gamma}'$ and the total increase in $d_{\gamma}'$.  {\bf Conclusion:} the total variation of exterior angle $ea(\caC)$ bounds the first-derivative of distance $d_{\gamma}$ to a chord as follows
\[
-ea(\caC)\le d_{\gamma}'(r)\le d_{\gamma}'(p)\,+\, \int_p^r (d_{\gamma}'')_+\,ds\le ea(\caC).
\]
Combined with the formula for the derivative of distance, it follows for $ea(\caC)$ small, that the angle between $\caC$ and the fibers of the projection to $\gamma$ are close to $\pi/2$.

We further consider the segment $\widehat{pq}$ of a chord, intersecting $\caC$ only at $p,q$ with: $r$ a distinct third point on $\caC$ between $p,q$, and $s$ the projection of $r$ to $\widehat{pq}$.  We now observe for $ea(\caC)$ small, that the projection of $r$ is an interior point of $\widehat{pq}$ and that the angle between $\widehat{pr}$ and $\widehat{rs}$ is close to $\pi/2$. The argument is as follows.  The union of $\widehat{pr}$ and the subarc $\caC_r$ of $\caC$ between $p,r$ forms a polygon $\mathcal P$ with $n+2$ vertices, for $n$ the number of vertices of $\caC_r$.  The sum $sa$ of interior angles of $\mathcal P$ at the vertices of $\caC_r$ satisfies $|sa-n\pi |\le ea(\mathcal C_r)$ and by Gauss-Bonnet the sum of all interior angles of $\mathcal P$ is bounded by $n\pi$.  It follows that the interior angle of $\mathcal P$ at $r$ is bounded by $ea(\caC_r)$.  From the paragraph above, the angle $\tilde\theta_r$ between the segment of $\caC_r$ containing  $r$ and $\widehat{rs}$ satisfies $|\cos \tilde\theta_r|\le ea(\caC_r)$.  Combining estimates, the angle $\theta_r$ between $\widehat{pr}$ and $\widehat{rs}$ is bounded as $|\theta_r-\pi/2|\le ea(\caC_r)+\arcsin ea(\caC_r)$.  {\bf Conclusion:} for $ea(\caC_r) \le 1/2$, then $|\theta_r-\pi/2| < \pi/2$ and consequently the chords $\widehat{pr}$ and  $\widehat{rq}$ are distinct from $\widehat{rs}$.  In particular $r$ projects to an interior point of $\widehat{pq}$ and $\widehat{rs}$ is orthogonal to $\widehat{pq}$.

The geodesic triangle $\Delta prs$ has a right-angle at $s$, and angle $\theta_r$ at $r$ close to $\pi/2$, with the difference bounded in terms of $ea(\caC_r)$.  The comparison triangle in the constant curvature $-\epsilon_{\caT}$ plane has  corresponding side lengths and angles at least as large, \cite[Chap. II.1]{BH}.  Since a sum of triangle angles is at most $\pi$, it follows for $ea(\caC_r)$ small, that the comparison triangle has two angles close to $\pi/2$ with differences bounded in terms of $ea(\caC_r)$.  We now recall that for triangles in constant negative curvature with two angles close to $\pi/2$ the length of the included side is correspondingly bounded.  In particular for a triangle $\Delta abc$ in hyperbolic geometry with angles 
$\alpha,\beta$ close to $\pi/2$ then the angle $\gamma$ is close to zero and the length of the opposite side satisfies
\[
\cosh c\,=\,\frac{\cos\alpha\cos\beta+\cos\gamma}{\sin\alpha\sin\beta}.
\]
The angle bounds combine to bound $\cosh c$ and $c$.  {\bf Conclusion:} for $ea(\caC_r)\le 1/2$, there is a universal bound for the distance between $\caC$ and any chord, with the bound small for $ea(\caC_r)$ small.

In general for a chord of a concatenation, we consider minimal segments connecting points of the concatenation and apply the above considerations to establish the following.
\begin{proposition}\label{distbd}
Let $\caC$ be a concatenation of geodesics with vertices in $\caT$ and total exterior angle variation $ea(\caC)\le 1/2$.  There is a uniform bound in terms of $ea(\caC)$  for the distance between a chord of $\caC$ and the corresponding subarc of $\caC$.  The bound for distance is small for $ea(\caC)$ small.
\end{proposition}

The counterpart result for concatenations with small exterior angle in the upper half plane is also valid.  

\subsection*{Main considerations}

We are ready to combine considerations.  We begin from Theorem \ref{bkrnd} with a bi sequence $\{\gamma_n\}$ of singular (cusp-to-cusp) unit-speed geodesics of $\mathcal T$ with dense image in the unit tangent bundle of $\mathcal T/MCG$.  We will apply elements of $MCG$ and perform {\bf twisting of geodesics} to obtain a concatenation $\caC$ with total exterior angle variation $ea(\caC)$ appropriately small.  We find that suitable chords of $\caC$ limit to a bi infinite geodesic $\caC_{\infty}$, strongly asymptotic to $\caC$. It follows that the geodesic $\caC_{\infty}$ is dense in the unit tangent bundle.  The result is presented in Theorem \ref{dense}.

We consider the sequence $\{\gamma_n\}$ and define the concatenation $\mathcal C$.  Select $\delta_1>0$ smaller than half the length of any singular geodesic and $\delta_2,\, \delta_1>2\delta_2>0$.  The cusps of $\caT$ lie in a single $MCG$ orbit.  Begin with the geodesic $\gamma_1$ and $MCG$ translate the geodesic $\gamma_2$ to arrange that the terminal point of $\gamma_1$ coincides with the initial point of the translated $\gamma_2$.   Perform a suitable twisting of $\gamma_2$ with distance parameter $\delta_2$ (see twisting of geodesics) to arrange for a concatenation with small exterior angles.  Proceed inductively, using the distance parameter $\delta_2$, to define an infinite concatenation $\mathcal C$ with: i) the sequence of exterior angles absolutely summable, and ii) all geodesic chords containing corresponding subarcs of $\caC$ within $\delta_2$ neighborhoods. The second property is realized by Proposition \ref{distbd}, for sufficiently small total exterior angle variation.  From properties of twisting of geodesics, since the twisting exponents tend to infinity along $\mathcal C$, the distance from $\mathcal C$ to the $MCG$-translates of $\{\gamma_n\}$ tends to zero and the pair have coinciding accumulation sets in $\mathcal T/MCG$.  The pair also have coinciding accumulation sets in the unit tangent bundle, a consequence of the $C^1$-approximation for twisting of geodesics. {\bf Conclusion:} the concatenation $\caC$ has dense image in the unit tangent bundle of $\mathcal T/MCG$

We note that geodesic chords of $\caC$ do not degenerate.  Write $B(\delta)$ for the $\delta$ neighborhood of a particular cusp.  We first show that the intersection $\caC_0=\caC\cap B(\delta_1)$ is a single connected segment.  By construction, $\caC_0$ is within $\delta_2$ of a segment of the concatenation of $MCG$-translates of $\{\gamma_n\}$, and by construction, $\caC_0$ is also within $\delta_2$ of the chord connecting its endpoints.  By convexity the entire chord is contained within $B(\delta_1)$.  The estimates provide that the segment of the $MCG$-translates of $\{\gamma_n\}$ is within $2\delta_2$ of the chord and within $\delta_1+2\delta_2$ of the cusp.  The closest distinct cusps are at distance strictly greater than $\delta_1+2\delta_2$.  It follows that the (maximal) segment within $B(\delta_1)$ of the $MCG$-translates of $\{\gamma_n\}$ consists of only two geodesics. Since $\caC$ and the concatenation of $MCG$-translates of $\{\gamma_n\}$ coincide outside of a $\delta_2$ neighborhood of any cusp, it follows that the concatenations agree on $B(\delta_1)-B(\delta_2)$.  The description of $\caC\cap B(\delta_1)$ follows.  The intersection of any chord of $\caC$ and $\partial B(\delta_1/2)$ is contained in a $\delta_2$ neighborhood of $\caC_0$.  The intersection of a $\delta_2$ neighborhood of $\caC_0$ and $\partial B(\delta_1/2)$ is compact, since $\delta_1>2\delta_2$. {\bf Conclusion:} chords of $\caC$ do not degenerate on Teichm\"{u}ller space.   

We construct a limit of chords. Let $\{p_n\}$ be a sequence of points along $\mathcal C$ with $p_n$ on $\gamma_n$.  By the compactness argument of \cite[Sec. 7]{Wlcomp} the sequence $\{\widehat{p_{-n}p_n}\}$ of chords has a subsequence which converges in a generalized sense (a priori the limit may be degenerate). From the above, the subsequence does not degenerate. {\bf Conclusion:} a limit  $\mathcal C_{\infty}$ of chords is a bi infinite geodesic on Teichm\"{u}ller space (in general a limit is described by a sequence of singular geodesics).   

We show that the distance between $\caC$ and $\mathcal C_{\infty}$ tends to zero at infinity.  The geodesic $\mathcal C_{\infty}$ is the base of a nearest point projection.  Consider the distance $F=d(\mathcal C_{\infty},*)$ to $\mathcal C_{\infty}$ as a function of arc length along $\mathcal C$.   The function is piecewise convex with total variation of discontinuities 
of $F'$ bounded by $ea(\caC)$.  Consider first that $F$ has an infinite number of zeros. The geodesic $\caC_{\infty}$ is approximated by chords, and so $F$ is approximated by the distance between $\caC$ and a sequence of chords. Thus between a pair of zeros of $F$, the function is bounded by Proposition \ref{distbd}.  Between zeros, $F$ is bounded by the sum of the absolute discontinuities of $F'$.  The absolute sum of discontinuities of $F'$ between zeros of $F$ tends to zero along $\caC$, and thus $F$ limits to zero at infinity.   Consider second that $F$ has only finitely many zeros.
Since $F$ is bounded, $0<F<\delta_2$, it follows that the total variation of $F'$ is bounded.  It follows that the terminal total variation of $F'$ tends to zero.  Equivalently stated, $F'$ tends to a limit, which in fact is zero since $F$ is bounded. 
By construction $\caC$ consists of segments inside $\delta_2$ neighborhoods of cusps, connected by geodesic segments (outside of cusp neighborhoods) of length at least $2\delta_1-2\delta_2>0$.  By the second variation of distance formula and strict negativity of curvature, there is a positive lower bound in terms of the magnitude of $F$ for the convexity of $F$ on the geodesic segments.  Since $F'$ has a limit along $\caC$, it follows that the convexity of $F$ tends to zero along $\caC$.  
It follows that $F$ limits to zero along $\caC$.  The curves $\mathcal C$ and $\mathcal C_{\infty}$ are strongly asymptotic, and in fact have coinciding accumulation sets in the unit tangent bundle, since close geodesic segments have tangent fields close for proper subsegments. We have the following.
\begin{theorem}
\label{dense}
A limit of chords between points tending respectively to positive and negative infinity along the concatenation $\caC$ is an infinite geodesic, dense in the unit tangent bundle of $\mathcal M=\caT/MCG$. 
\end{theorem}

The corresponding considerations for the upper half plane are also valid.  There are two basic  matters to address to apply the approach for higher dimensional moduli spaces.  The first is that the strict convexity of $CAT(-\epsilon)$ geometry is replaced by the convexity of $CAT(0)$ geometry. For the latter geometry, there is no bound for distance only in terms of exterior angle; an alternative to Proposition \ref{distbd} is needed.  The second is the non trivial  condition for geodesics emanating from a cusp to be approximated by twisting geodesics, see the above comments before {\bf twisting geodesics}. Addressing the matters is an open question.

\section{Recurrence, volume growth, horseshoes, and closed geodesics}
The main result of the section is that the moduli space $\mathcal M$ of once punctured tori contains infinitely many WP closed geodesics, whose number grows exponentially in length.  Although by a comparison of metrics the result follows from the corresponding result for the Teichm\"{u}ller metric, we present an argument not involving the Teichm\"{u}ller geometry.  

We begin by studying the recurrence of the geodesic flow $\phi_T \colon S\mathcal M \to S\mathcal M$ on the unit tangent bundle.  A vector $v \in S\mathcal M$ is {\em recurrent} if for every neighborhood $U$ of $v$ there exists arbitrarily large $T > 0$ such that $\phi_T v \cap U \ne \emptyset$. 

\begin{lemma} Almost every vector in $S\mathcal M$ is recurrent.
\end{lemma}

\prf Wolpert \cite{Wlpi} showed that WP $\text{Area}(\mathcal M) = \pi^2/12 < \infty$ and thus $\text{Vol}(S\mathcal M) <\infty$. (Mirzakhani's general recursion \cite{Mirvol} determines all values $\text{Vol}(\mathcal M_{g, n})$.)  The GF  preserves the Liouville measure and from the Poincar\'e recurrence theorem it follows that for any open set $U \subset S\mathcal M$ we can find a $T>0$ with $\phi_T U \cap U = \emptyset$.  The result follows.

The topological entropy $h$ of a complete geodesic flow on a {\em compact surface} is a measure of the exponential rate of divergence of nearby geodesics.  The entropy $h$ provides one measure of the complexity of the global geodesic structure. Another measure of complexity is the exponential growth rate of closed geodesics $h_{CG}$, defined as the exponent coefficient for the exponential growth rate in $T$ of the number of closed geodesics with length $\leq T$.  A third measure is the volume entropy $h_{V, x}$ defined as the limit $\lim_{r \to \infty} (1/r) \log Area(B(x,r)),$ where $B(x, r) $ denotes the ball of radius $r$ around the point $x$ in $\mathcal T$. For a compact surface with everywhere negative curvature, the geodesic flow is Anosov \cite{Anosov} and the three entropies coincide, i.e. $h = h_{CG} = h_{V, x}$ for every $x$ \cite{Bowen1, Mann}. We will see that the topological entropy for the GF is positive while the volume entropy is infinite.  

For a compact surface a geodesic flow with positive topological entropy contains a horseshoe \cite{Katok}. The same statement holds for the GF restricted to a compact invariant subset.  Although a horseshoe in $S\mathcal M$ could have zero volume, the dynamics of the GF on the horseshoe would be chaotic. In particular the horseshoe would contain infinitely many closed geodesics, whose number would grow exponentially in length. We show the following.

\begin{proposition} 
\label{expgr}
The moduli space $\mathcal M$ contains infinitely many closed geodesics, whose number grows exponentially in length. 
\end{proposition}

Since the flow $\phi_T \colon S \mathcal M \to S\mathcal M$ is {\em not} complete, there is no natural definition of topological entropy.  We will describe a compact, invariant subset of $S \mathcal M$ on which the GF has positive topological entropy.  This will imply $h_{CG} > 0$ and the desired conclusion. \\

\prf For $\epsilon > 0$, let $S\mathcal M(\epsilon)$ be the collection of unit tangent vectors whose corresponding geodesics never enter the $\epsilon$-neighborhood of the cusp.  The set $S\mathcal M(\epsilon)$ is compact and $\phi_T \colon S\mathcal M(\epsilon) \to S\mathcal M(\epsilon)$ is complete.  Since the tangent vectors in  $S\mathcal M(\epsilon)$ only visit regions of uniformly bounded negative curvature, $S\mathcal M(\epsilon)$ is a non-maximal hyperbolic set for the GF 
\cite{Brin}. 

In general a Markov partition provides a symbolic model for a smooth flow  as a special flow over a subshift of finite type \cite{Brin}.  Bowen showed that every locally maximal, hyperbolic set has a Markov partition \cite{Bowen}. More recently, Fisher \cite{TF} showed that every hyperbolic set has arbitrarily small enlargements on which the flow has a Markov partition.\footnote{Fisher established the result for diffeomorphisms, but his proof immediately extends to sections of flows.} For the WP this implies that for each sufficiently small $\epsilon>0$, there is a hyperbolic set compact enlargement $V(\epsilon) \supset S \mathcal M(\epsilon)$, which stays away from the cusp, and for which  
$\phi_T \colon V(\epsilon) \to V(\epsilon)$ has a Markov partition.

The flow $\phi_T \colon V(\epsilon)\to V(\epsilon)$  has positive topological entropy provided that the associated subshift of finite type has positive topological entropy \cite{Abram}. The subshift of finite type is defined by an adjacency matrix: a square matrix having non-negative entries. The structure theorem for adjacency matrices \cite{LM} yields that either the subshift has only finitely many periodic points or has positive topological entropy (and thus infinitely many periodic points, whose number grows exponentially in the period.)  It follows that in  $V(\epsilon)$, either the GF has only finitely many closed geodesics or contains infinitely many closed geodesics, whose number grows exponentially in length. 

The following geometric lemma provides that some $V(\epsilon)$ contains infinitely many simple closed geodesics.   The lemma completes the proof of Proposition \ref{expgr}.  
\begin{lemma}
There is a compact subset of $\mathcal M$ containing an infinite number of closed geodesics.
\end{lemma}
\prf  Choose a torsion free subgroup $\Gamma$ of $MCG$ with $\mathcal T/\Gamma$ having a non trivial, non peripheral simple closed curve. From the classification of surfaces the surface $\mathcal T/\Gamma$ has infinitely many distinct non peripheral, free homotopy classes of simple closed curves.  From the classification of elements of $MCG$ and the WP axis theorem \cite[Theorem 25]{Wlcomp} each non trivial, non peripheral free homotopy class contains a closed geodesic. From the minimal intersection property for negative curvature each representing geodesic is simple. 

We observe that simple WP geodesics are uniformly bounded away from cusps and so are contained in a compact subset of 
$\mathcal T/\Gamma$. The lift of a simple geodesic to $\mathcal T$ is disjoint from its Dehn twist translates.  In FN coordinates $(\ell,\tau)$ and for $n_0$ the smallest exponent with $T^{n_0}\in \Gamma$, a simple geodesic has a lift intersecting a neighborhood of the cusp $\{\ell=0\}$ and contained entirely in the sector $\{0<\tau<3n_0\ell\}$.  The lift intersects $\{\ell=1\}$ in the compact set $\{\ell=1,\,0<\tau<3n_0 \}$.  For a $CAT(-\epsilon)$ geometry geodesic segments depend continuously on endpoints.  The geodesics intersecting the compact set are bounded away from $\{\ell=0\}$, as desired.   The proof is complete.

Finally, we note that the volume entropy  $h_{V, x}$ is infinite for every point.  This holds since the area of a sufficiently large WP ball in $\mathcal T$ is infinite.  A sufficiently large ball contains a cusp point and the WP area of a neighborhood in $\Tbar$ of a cusp point is infinite. In particular the WP area form is $\omega=\frac12\,d\ell \wedge d\tau$ and a small metric neighborhood of a cusp point is parameterized in FN coordinates by  
$0 < \ell < \epsilon,\, -\infty < \tau < \infty$.


\end{document}